\documentclass[a4paper,runningheads]{llncs}
\usepackage[utf8]{inputenc}
\usepackage{booktabs} %
\usepackage[dvipsnames]{xcolor}
\usepackage{relsize}
\usepackage{amsmath,amssymb,amsfonts} %
\renewcommand{\vec}[1]{\mathbf{#1}}
\usepackage{algorithm}
\usepackage[noend]{algorithmic}

\algsetup{linenosize=\tiny}

\usepackage{pdflscape}
\usepackage{amsmath}

\usepackage{hyperref}
\hypersetup{colorlinks, allcolors=black, urlcolor=blue}

\usepackage[square,numbers,sort]{natbib}
\makeatletter
\renewcommand\bibsection%
{\section*{\refname\@mkboth{\MakeUppercase{\refname}}{\MakeUppercase{\refname}}}}
\renewcommand\@biblabel[1]{#1.}

\makeatother

\usepackage{subfigure}
\usepackage{graphicx}
\usepackage{array}

\DeclareMathOperator*{\argmin}{arg\,min}

\newcommand{\assign}{\ensuremath{:=}}
\newcommand{\procedure}[1]{\textsf{#1}}

\newcommand{\mcol}{\multicolumn}

\usepackage{orcidlink}
\renewcommand{\orcidID}[1]{\raisebox{1.1ex}{\orcidlink{#1}}}
\let\svthefootnote\thefootnote
\newcommand\blankfootnote[1]{%
  \let\thefootnote\relax\footnotetext{#1}%
  \let\thefootnote\svthefootnote%
}

\begin{document}
\blankfootnote{\smaller Please cite as: M.\@ López-Ibáñez, F.\@ Chicano, R.\@ Gil-Merino. The Asteroid Routing Problem: A Benchmark for Expensive Black-Box Permutation Optimization. In J. L. Jiménez Laredo et al., editors, \emph{Applications of Evolutionary Computation}, volume 13224 of Lecture Notes in Computer Science. Springer Nature, Switzerland, 2022. \url{https://doi.org/10.1007/978-3-031-02462-7_9}}

\title{The Asteroid Routing Problem: A Benchmark for Expensive Black-Box Permutation Optimization}
\titlerunning{The Asteroid Routing Problem}
\author{Manuel~L{\'o}pez-Ib{\'a}{\~n}ez\orcidID{0000-0001-9974-1295} \and %
  Francisco~Chicano\orcidID{0000-0003-1259-2990} \and %
  Rodrigo Gil-Merino}
\authorrunning{M. L{\'o}pez-Ib{\'a}{\~n}ez, F. Chicano and R. Gil-Merino}
\institute{ITIS Software, Universidad de Málaga, Spain \\
\email{manuel.lopez-ibanez@uma.es, chicano@uma.es, gilmerino@uma.es}}

\maketitle              %
\begin{abstract}
  Inspired by the recent 11th Global Trajectory Optimisation Competition, this
  paper presents the asteroid routing problem (ARP) as a realistic benchmark of
  algorithms for expensive bound-constrained black-box optimization in
  permutation space. Given a set of asteroids' orbits and a departure epoch,
  the goal of the ARP is to find the optimal sequence for visiting the
  asteroids, starting from Earth's orbit, in order to minimize both the cost,
  measured as the sum of the magnitude of velocity changes required to complete
  the trip, and the time, measured as the time elapsed from the departure epoch
  until visiting the last asteroid. We provide open-source code for generating
  instances of arbitrary sizes and evaluating solutions to the problem.  As a
  preliminary analysis, we compare the results of two methods for expensive
  black-box optimization in permutation spaces, namely, Combinatorial Efficient Global Optimization (CEGO), a Bayesian
  optimizer based on Gaussian processes, and Unbalanced Mallows Model (UMM), an estimation-of-distribution
  algorithm based on probabilistic Mallows models. We investigate the best
  permutation representation for each algorithm, either rank-based or
  order-based. Moreover, we analyze the effect of providing a good initial
  solution, generated by a greedy nearest neighbor heuristic, on the
  performance of the algorithms. The results suggest directions for
  improvements in the algorithms being compared.
  
\end{abstract}

\section{Introduction}

Several space programs from different countries focus on small bodies orbiting around the Sun and the Earth. Apart from the search for answers about the origin and evolution of the Solar System, there are other practical reasons that space agencies consider. One of these reasons is related to the present decrease of mineral resources on Earth. These resources, heavily used in technological devices like mobile phones and computers, include silicon, quartz, boronite and others, and are hard to mine for: the density of these minerals is usually high, and they fell to inner strata during Earth formation, which translates into an extra difficulty when mining for these materials.
As an alternative to Earth mining for techno-minerals, space agencies have proposed asteroids and near Earth objects (NEOs) \cite{esahera,esaasteroid} as near-term mining targets. These objects contain varying amounts of rare earth metals, minerals and water at surface levels, consequently much easier to mine for than on Earth. A clear drawback to these plans is the transport from and to the asteroids, which will imply high consumption of energy.
A solution to reduce the energy consumption is to build human settlements on large asteroids and/or on orbit, together with space solar power stations that could fuel spacecrafts and those settlements. This idea firstly appeared on a science-fiction novel and later it was popularized as a ``gedankenexperiment" by Dyson in the 60's. Nowadays, the idea of a swarm of space solar power stations is called ``the Dyson sphere".
In its last edition, the \emph{11$^\text{th}$ Global Trajectory Optimisation Competition} (GTOC11, \url{https://gtoc11.nudt.edu.cn}) proposed to computationally implement a Dyson ring, which is a set of space stations (12 in this case) orbiting the Sun in a circular orbit, to mine close to Earth passing-by asteroids. The Dyson ring construction was proposed as an optimisation problem in the context of orbital mechanics. The first part of the optimization problem consists in finding a plan (trajectory) to visit a selected subset of asteroids with a spacecraft launched from Earth. This spacecraft does not return to Earth. The second part of the problem consists in activating a device, called \emph{asteroid transfer device} or ATD, located in the visited asteroids to guide them to the space stations in the Dyson ring. Once the Dyson ring is formed, regular visits to the space stations in the Dyson ring are supposed to bring the minerals in the asteroids to Earth. We are interested in this paper in solving the first optimization problem, for which we provide a formulation as a challenging expensive black-box permutation problem that we call \emph{Asteroid Routing Problem} (ARP) that may be used to compare and benchmark optimization approaches.

In recent years, there has been growing interest in tackling expensive black-box
permutation problems~\citep{ZaeStoFriFisNauBar2014,PerLopStu2015si,BarZae2017model,IruLop2021gecco}. In permutation problems, the decision space is the space
of permutations of a given length $n$, usually denoted with $S_n$. Black-box optimization, as considered in this paper, tackles an optimization problem without requiring an explicit
analytical model of its objective function and constraints. Finally, in
expensive optimization, each function evaluation is often the result of a
costly simulation or physical experiment, thus the total budget of evaluations
allowed before returning the best solution found is severely restricted,
typically no more than 1000~\citep{KnoCorRey09emo,PerLopStu2015si}. The
combination of these three characteristics gives rise to a challenging family
of optimization problems. Typical approaches include estimation{-}of{-}distribution
algorithms~\citep{IruLop2021gecco} and combinatorial Bayesian
optimizers~\citep{ZaeStoFriFisNauBar2014} or combinations of both~\citep{PerLopStu2015si}.

A potential criticism of recent works is that the benchmark problems considered ---e.g.,
 the travelling salesman problem, the linear ordering problem, etc.--- are
taken from the classical combinatorial optimization literature and they are
neither truly expensive nor black-box.
The ARP proposed here is much more expensive to evaluate than classical
combinatorial optimization problems, even for short permutations. Moreover,
there is no explicit analytical model of the problem from the routing
perspective, thus classical techniques such as fast neighborhood evaluation or
gradient calculation are not possible. In our formal definition of the ARP,
we aim to preserve the most challenging features of the problem, without
requiring expert knowledge of the astrophysical details. Benchmark
generation and evaluation code is made publicly available.  We also provide
preliminary results using state-of-the-art optimization methods.  Our hope is that
this work will raise interest in the field of expensive black-box combinatorial
optimisation and its applications to trajectory optimisation.

This paper is structured as follows. Section~\ref{sec:background} gives a brief
introduction to basic astrodynamical concepts required to understand the
problem as an optimization benchmark. We do not explain in detail the
underlying calculations because (1) these are standard formulas available in
many physics textbooks, (2) they are implemented by multiple open-source
software libraries and (3) they are not directly available to any of the
optimization algorithms considered here. The optimization model of
the proposed Asteroid Routing Problem is given in Section~\ref{sec:problem}. The optimization algorithms evaluated in this work are described in Section~\ref{sec:algos}. We provide an experimental analysis of these algorithms on the ARP in Section~\ref{sec:analysis} and conclude in Section~\ref{sec:conclusion} with a summary of findings and suggestions for further research motivated by these findings as well as possible extensions of the ARP benchmark. 

\section{Background}\label{sec:background}

In this section we will briefly revise the basic theory behind the astrodynamics of the problem. 
We will start describing the two-body problem and its solution leading to the Keplerian orbits. Then, we will describe one type of maneuver used by the spacecrafts to change their trajectory in space and we will end this section presenting the Lambert problem to determine the impulses required to reach a point in space and time from any other point.

\subsection{Two-body Problem}

The dynamics of two bodies in space subject to the gravitational force is determined by Newton's second law of dynamics and gravitation law. Combined, they lead to the following vector equation:
\begin{equation}
\label{eqn:two-body}
\frac{d^2\vec{r}}{dt^2} = - \frac{\mu}{r^3} \vec{r},
\end{equation}
where $\vec{r}$ is the position vector of the second object respect to the first one, $r = |\vec{r}|$ is the module of the position vector, $t$ is time, and $\mu$ is usually called \emph{gravitational parameter}, which for our solar system  is typically assumed to be a constant due to the mass of the Sun being much larger than the mass of any other object.

Despite the dynamics of the two bodies happening in 3D space, the conservation of the angular momentum forces the movement to happen in a plane and the analytical solutions of the problem are conic curves (also called Keplerian orbits when they refer to the two-body problem): ellipse (and circumference as a particular case), hyperbola and parabola\footnote{The trajectory of the second object could also be linear, but this is not of practical interest, because this would end with the object destroyed in the Sun.}. Both, the hyperbola and parabola are open curves, which means that the object following one of them will come from outside the solar system and after approaching the Sun will disappear outside the solar system without coming back again. Solar system objects, like planets and asteroids, follow an elliptic (or circular) orbit around the Sun. 
\begin{figure}[!htb]
	\centering
	\includegraphics[scale=0.30]{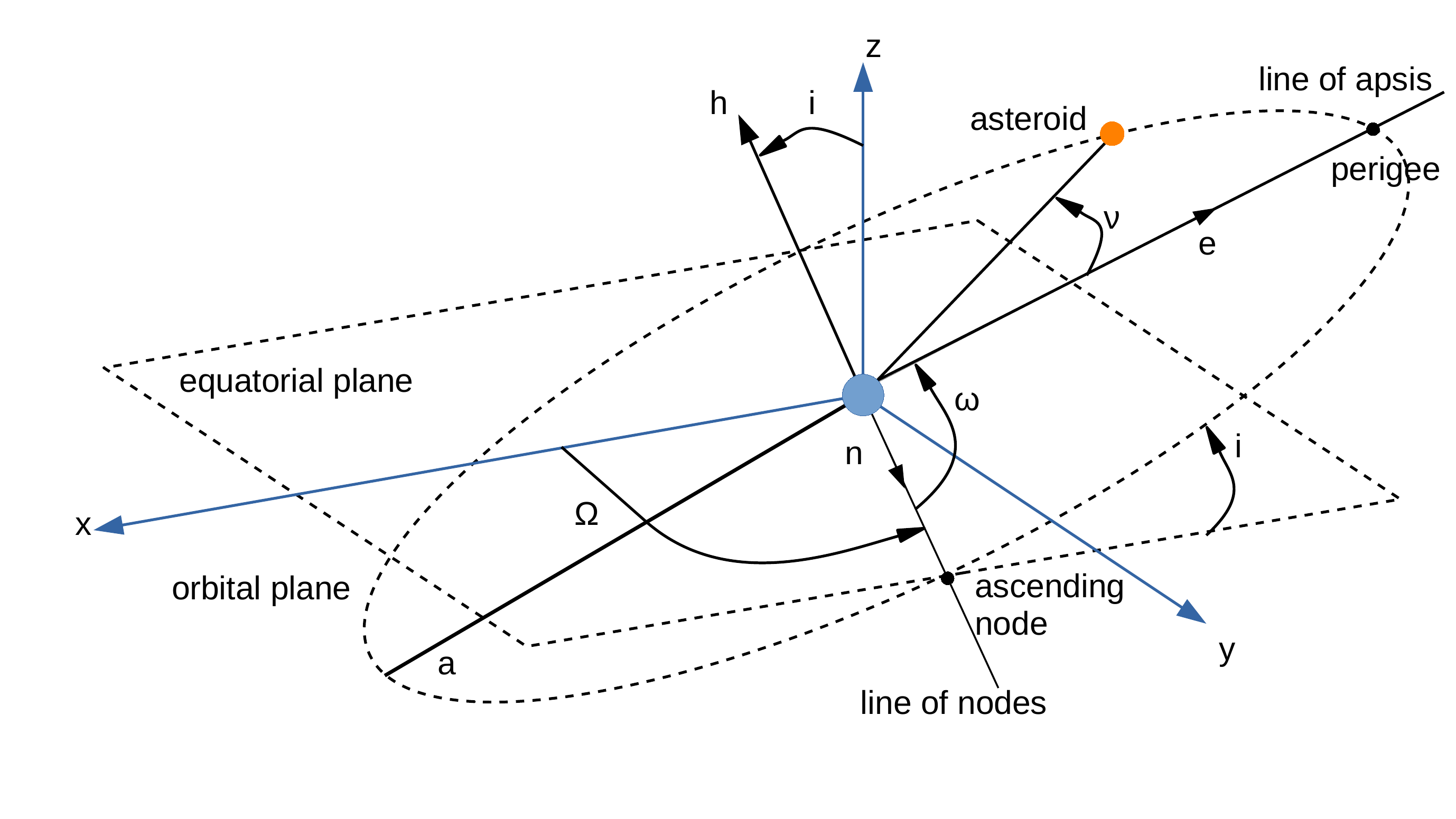}
	\caption{The five parameters to describe an elliptic orbit.}
	\label{fig:orbit_param}
\end{figure}

When we fix a reference frame in space, elliptic orbits can be characterized by five \emph{orbital elements}: ${a}$ (semi-major axis), $e$ (eccentricity), $i$ (inclination), $\omega$ (argument of periapsis\footnote{The periapsis is the nearest point of an elliptic orbit to the object in the focus. When the object is the Sun, the periapsis is also called perihelion.}) and $\Omega$ (longitude of the ascending node) (see Figure~\ref{fig:orbit_param}). The position of the object in the orbit is characterized by a sixth parameter which changes with time: the true anomaly, $\nu$. The term \emph{epoch} is used in this context to refer to a specific \emph{date} and \emph{time}. The movement of an object, i.e., its position and velocity vectors, in an elliptical orbit is completely determined by the six orbital elements described above in a given epoch.

When there are more than two bodies, as it happens in our solar system, the movement equations cannot be solved analytically, except for a few particular cases. In that non-analytical case, the solution of the two-body problem is still useful. The presence of additional bodies can be considered in most of the cases a minor perturbation in the Hamiltonian describing the two-body problem and, as a consequence, the orbital elements are slowly changing with time \cite[Chapter 11]{goldstein:mechanics}. In our formulation of the problem, however, we do not consider these perturbations, and we assume that all objects follow a Keplerian orbit.

\subsection{Maneuvers in Space}

Artificial spacecrafts also follow a Keplerian orbit around the Sun when they are far from the gravitational influence of the Earth or any other solar system object. They usually have an engine that allows them to change the orbit using operations called \emph{maneuvers}. There are two main kinds of maneuvers: \emph{continuous} and \emph{impulsive}. In a continuous maneuver the engine of the spacecraft is on during a long period of time (e.g., days) and generates a force whose direction and magnitude can be controlled. We will only consider impulsive maneuvers here.

Impulsive maneuvers cover a wide range of different space maneuvers, from small orbital adjustments when two spacecrafts encounter each other, to large orbital corrections between a low-Earth orbit and an interplanetary transit. To carry out an impulsive maneuver, usually small engine firings take place, changing only the direction and magnitude of the velocity vector instantaneously. This means that the position vector does not change during the impulse, which, although an idealisation, it works in most situations. The advantage of treating the impulsive maneuver in this way is that it avoids solving the equations of motion including the thrusts.
Under these assumptions, the trajectory of the spacecraft is still a Keplerian orbit where the orbital elements change instantaneously at the point where the impulsive maneuver happens.

\subsection{Lambert Problem}

At the heart of Astrodynamics, a key problem is to determine the Keplerian orbit connecting two positions in a given period of time $\Delta t$. This is the Lambert's problem, also known as the Gauss' problem~\citep{Izz2014lambert}.

We can use the solution of Lambert's problem to determine the impulsive maneuvers required for a spacecraft to do a rendezvous with an asteroid. The transit orbit and impulsive maneuvers can be completely characterized by two times: the epoch at which the transit starts $\tau$ and the transit time $t$. These two times determine the two points $P_1$ and $P_2$ in space we need to join. $P_1$ is the point at which the spacecraft is at $\tau$. $P_2$ is the point at which the asteroid is at $\tau+t$. Lambert's problem is solved to go from $P_1$ to $P_2$ in time $t$. Two impulsive maneuvers are used to move from $P_1$ to $P_2$. A first impulsive maneuver $\Delta \vec{v}_1$ at $P_1$ inserts the spacecraft in the transit orbit determined by solving Lambert's problem. A second impulsive maneuver $\Delta \vec{v}_2$ at $P_2$ inserts the spacecraft in the same orbit as the asteroid (rendezvous). Let's denote with $s$ the orbit of the spacecraft (six orbital parameters in an epoch) and with $a$ the orbit of the asteroid\footnote{We previously used $a$ to mean the semi-major axis of an orbit, from now on we will use it to denote the complete orbit of an asteroid.}. Then, $\Delta \vec{v}_1$ and $\Delta \vec{v}_2$ are determined by $s$, $a$, $\tau$ and $t$:
\begin{equation}
\label{eqn:lambert}
    (\Delta \vec{v}_1, \Delta \vec{v}_2) = \text{Lambert} (s, a, \tau, t)
\end{equation}

\section{Asteroid Routing Problem}
\label{sec:problem}

The \emph{Asteroid Routing Problem} consists in finding a route followed by a spacecraft launched from Earth to visit a given set of $n$ asteroids $A=\{a_1, a_2, \ldots, a_n\}$ that minimizes both the sum of velocity impulses required by the route (related to consumption) and the total time required to visit all of them. A solution to this problem is a pair $(\pi, \vec{t})$, where $\pi \in S_n$ is a permutation representing the order in which the asteroids are visited and $\vec{t} \in \mathbb{R}^{2n}_{\geq 0}$ is a vector of $2n$ non-negative real numbers representing parking and transit times to reach each asteroid. 

We will consider a start epoch $\tau_0$, determined by the instance, when the spacecraft is on Earth. The spacecraft remains on Earth during a time $t_1$ after $\tau_0$. Then, it is launched at epoch $\tau_0 + t_1$ to reach the first asteroid in the sequence, $a_{\pi(1)}$, using a transit orbit computed by solving Lambert's problem between the Earth and the asteroid. In the launch, the spacecraft experiments an impulse $\Delta\vec{v}_1$ that inserts the spacecraft in a transit orbit to reach asteroid $a_{\pi(1)}$ in time $t_2$ after the launch epoch. Once the asteroid is reached, a new impulse $\Delta\vec{v}_2$ is used to make the spacecraft to follow the orbit of asteroid $a_{\pi(1)}$. 

In order to visit the $i$-th asteroid, $a_{\pi(i)}$, for $i>1$, the spacecraft remains during time $t_{2i-1}$ in the orbit of asteroid $a_{\pi(i-1)}$ and receives an impulse $\Delta \vec{v}_{2i-1}$ to reach asteroid $a_{\pi(i)}$ in time $t_{2i}$. Then, a new impulse $\Delta \vec{v}_{2i}$ changes the orbit of the spacecraft to coincide with the one of $a_{\pi(i)}$.

For the transit between asteroids $a_{\pi(i-1)}$ and $a_{\pi(i)}$, the impulses $\Delta\vec{v}_{2i-1}$ and $\Delta\vec{v}_{2i}$ are computed using Equation~\eqref{eqn:lambert}.
All times at odd positions of the time vector, $t_{2i-1}$, are parking times, because the spacecraft is waiting in an asteroid (or Earth's) orbit, while the times at even positions of the time vector, $t_{2i}$ are transit times, because the spacecraft is traveling through a transit orbit to reach the next asteroid.

The two objectives of the problem to minimize are the sum of all the velocity
impulses, which is related to energy consumption, and the sum of all parking
and transit times, which gives the elapsed time from launch to arriving at the
last asteroid. More formally, the two objectives are:
\begin{equation}\label{eq:objs}
  \Delta v = \sum_{i=1}^{2n} \left\lvert \Delta \vec{v}_i \right\rvert\quad\text{and}\quad
  T = \sum_{i=1}^{2n} t_i\enspace,
\end{equation}
where $\Delta \vec{v}_i$ is computed as follows:
\begin{equation}
  (\Delta \vec{v}_{2i-1}, \Delta \vec{v}_{2i}) = \text{Lambert}\left(a_{\pi(i-1)}, a_{\pi(i)}, \tau_{i-1} , t_{2i}\right)\enspace,
\end{equation}
where $\tau_{i-1} = \tau_0 +\sum_{j=1}^{2i-1} t_j$ is the launch epoch from
orbit $a_{\pi(i-1)}$ and we consider that $a_{\pi(0)}$ is Earth's orbit.

An instance of the problem is completely determined by the orbital parameters at a given epoch of Earth and all asteroids,
the starting epoch $\tau_0$
and the gravitational parameter $\mu$.

In this paper, we will focus on a scalarized variant of the above bi-objective problem that simply aggregates the two objectives in Eq.~\eqref{eq:objs} 
as follows:
\begin{equation}\label{eq:objfun}
  f(\pi,\vec{t}) = \Delta v + \frac{2~\text{km/s}~}{30~\text{days}~}\cdot T,
\end{equation}
where the constant in front of $T$ is based on preliminary experiments.

\section{Optimization Algorithms}\label{sec:algos}

In order to use the ARP as a benchmark problem for expensive black-box
combinatorial optimization, we tackle the problem in a hierarchical manner. The
\emph{inner} optimization problem decides the values of the vector of times
given a permutation of the asteroids, while the \emph{outer} problem aims to
find an optimal permutation of the asteroids. Since our focus here is the outer
problem, we always tackle the inner problem using the deterministic SLSQP (Sequential Least Squares Programming)
algorithm, thus ensuring that we obtain the same objective function
value given the same permutation. For the outer problem, we evaluate and
compare two algorithms for expensive black-box optimization problems, namely,
Unbalanced Mallows Model (UMM)~\citep{IruLop2021gecco} and Combinatorial Efficient
Global Optimization (CEGO)~\citep{ZaeStoFriFisNauBar2014}. In addition, we study the effect that a good initial solution has
on the behavior of the black-box algorithms. This initial solution is found by
means of a gray-box greedy nearest neighbor heuristic.  All these algorithms
are described below.

\subsection{Sequential Least Squares Programming (SLSQP)}

The inner optimization of the vector of times is done using the implementation
of SLSQP provided by
SciPy~\citep{Scipy2020natmet} without an explicit gradient, which is estimated
by the SLSQP algorithm using 2-point finite differences with an absolute step
size of 1.49e-08. Each transfer is optimized independently, that is, SLSQP
solves a sequence of $n$ problems with two numerical decision variables
$(x_0, x_1)$, which correspond respectively to the parking and transit times to
visit a particular asteroid. Parking times are bounded by $[0,730]$, while
transit times are bounded by $[1,730]$, both measured in days. 

\subsection{Greedy Nearest Neighbor Heuristic}

A reasonably good permutation of the asteroids can be generated by iteratively
visiting the nearest asteroid, in Euclidean distance, to the last-visited
asteroid, after calculating the positions of all unvisited asteroids at the
arrival time at the last-visited asteroid. This heuristic is presented in
Algorithm~\ref{alg:heuristic}. The algorithm keeps track of the orbit of the
spacecraft $s$ and its arrival time (epoch) $\tau$ at the last-visited
asteroid, which are initialized to the Earth's orbit $a_{\pi(0)}$ and the
initial epoch $\tau_0$, and a set of unvisited asteroids $U$. At each iteration
$i$, an asteroid $a_{\pi(i)}$ is chosen from $U$ according to the minimum value
of $d(s, a, \tau)$, which calculates the Euclidean distance between the
Cartesian positions of the orbits of the spacecraft $s$ and asteroid $a$ when both are
considered at time $\tau$.  After the asteroid $a_{\pi(i)}$ is chosen, we solve the inner
problem---using SLSQP as described above---to find the parking time $t_{2i-1}$
at asteroid $a_{\pi(i-1)}$ and transit time $t_{2i}$ to asteroid
$a_{\pi(i)}$. The algorithm finishes when all asteroids are visited, which
results in a permutation $\pi$ of the asteroids and the corresponding vector
$\vec{t}$ of parking and transit times.

\begin{algorithm}[!ht]
\caption{Greedy Nearest Neighbor Heuristic.}\label{alg:heuristic}
\begin{algorithmic}[1]
\STATE $s \assign a_{\pi(0)}$ \COMMENT{Earth's orbit} 
\STATE $\tau \assign\tau_0$ \COMMENT{Epoch of the spacecraft}
\STATE $U \assign \{1,\dotsc,n\}$
\FOR {$i \assign 1$ \textbf{to} $n - 1$}
\STATE $\pi(i) \assign \argmin_{j \in U} d(s, a_j, \tau)$\COMMENT{Euclidean distance between orbits at epoch $\tau$}
\STATE $(t_{2i-1}, t_{2i}) \assign \procedure{SLSQP}(s$, $a_{\pi(i)})$\COMMENT{Solve inner problem}
\STATE $\tau \assign \tau + t_{2i-1} + t_{2i} $
\STATE $U \assign U \setminus \{\pi(i)\}$
\STATE $s \assign a_{\pi(i)}$
\ENDFOR
\RETURN $(\pi, \vec{t})$
\end{algorithmic}
\end{algorithm}

\vspace{-12pt}

\subsection{Unbalanced Mallows Model (UMM)}

The UMM algorithm~\citep{IruLop2021gecco} is an
estimation-of-distribution algorithm based on the Mallows model and unbalanced
Borda learning.  The Mallows model for permutation spaces is defined by a
reference permutation $\sigma_0$ and a dispersion parameter $\theta$, which are
analogous to the mean and variance of the Gaussian distribution.
At each step of UMM, the permutations evaluated so far together with their
objective values are used to calculate a \emph{weighted mean} permutation
$\hat{\sigma}_0$.  From this $\hat{\sigma}_0$ and a value of the dispersion
parameter $\theta$, a new permutation is sampled at each iteration from the
Mallows model. The dispersion parameter $\theta$ is set by UMM in such a way
that, at the first iteration, the expected Kendall's-$\tau$ distance  of newly sampled
permutations from the reference permutation is at half of the expected distance
of an uniform sample. In subsequent iterations, the parameter $\theta$ is
adjusted such that the expected Kendall's-$\tau$ distance linearly decreases until the expected
distance is 1 at the last iteration.

Empirical results of UMM on expensive black-box variants of the linear ordering
problem and the permutation flowshop scheduling problem showed
results~\citep{IruLop2021gecco} comparable to the Combinatorial Efficient
Global Optimization (CEGO), with the additional benefit of being
computationally much less expensive than CEGO.

\subsection{Combinatorial Efficient Global Optimization (CEGO)}

CEGO extends classical Bayesian
optimizers based on Gaussian processes to unconstrained black-box combinatorial
optimization problems~\citep{ZaeStoFriFisNauBar2014}. A Gaussian process may be
  used as a surrogate of the continuous landscape of an expensive optimization
  problem. In order to model the landscape of a permutation space, CEGO uses a
distance measure for permutation spaces, such as the Kendall's-$\tau$ distance,
to \emph{interpolate} between solutions.
CEGO uses a genetic algorithm (GA) to explore the landscape of the surrogate
model. The GA optimizes the expected improvement criterion, which
takes into account both the expected mean and variance of a solution.  The
solution returned by the GA is then evaluated on the true objective function and
this information is used to update the Gaussian process surrogate model.

\section{Experimental Study}\label{sec:analysis}

In this section we describe the experiments performed and the results obtained. Our goal is to answer two research questions:

\begin{description}
 \item[RQ1:] How do the algorithms perform on the problem in a black-box setting (without any problem knowledge)?
 \item[RQ2:] How do they perform when some a priori information about the problem,
   in the form of a good initial solution, is introduced?
 \end{description}

In Subsection~\ref{subsec:methodology} the methodology used to answer both research questions is presented. Then, the experimental results to answer RQ1 and RQ2 are presented and described in Subsections~\ref{subsec:black-box} and~\ref{subsec:heuristic}, respectively.

\subsection{Experimental Methodology}
\label{subsec:methodology}

We prepared an instance generator that, given an instance size $n$ and seed for
the random number generator, randomly selects $n$ asteroids from the 83\,453
asteroids provided by the GTOC11 competition to create an instance of the
ARP. In this manner, we generated two instances (seeds 42 and 73) of each size
$n=\{10,15,20,25,30\}$. In the remainder of the paper, instances are named $n$\_\textit{seed}, that is, instance \texttt{10\_73} is an instance of size 10 generated with seed 73.

In our first set of experiments, designed to answer RQ1, both UMM and CEGO
start the search from a set of 10 initial solutions generated by a
\emph{max-min-distance sequential design}, that is, solutions are added iteratively to
a set of existing solutions by maximizing the minimum distance to solutions
already in the set. The first solution in the set is generated uniformly at
random. Generating the max-min-distance design only takes into account the
distance between permutations and, thus, it does not require any evaluation of
the objective function. These 10 solutions are then evaluated on the objective function
function (Eq.~\ref{eq:objfun}) and they become the initial population of either CEGO or UMM.
In order to answer RQ2, we initialize the max-min-distance design with the
heuristic solution returned by the greedy nearest neighbor heuristic described
above.  

Both UMM and CEGO stop after evaluating 400 permutations using the objective
function (Eq.~\eqref{eq:objfun}).  Each evaluation requires solving the inner problem by optimizing
the vector of times $\vec{t}$ using a sequence of $n$ runs of
SLQSP. Each run $i$ of SLSQP optimizes the pair of parking and transit times
$(t_{2i-1}, t_{2i})$, starting from the initial solution $(0,30)$ and
performing a maximum of 1000 iterations.

For each algorithm, we consider
two representations, namely, whether the permutations are represented as the
\emph{order} or as the \emph{ranking} of the asteroids in the visiting
sequence. The \emph{order}-based representation is the one used in the description of the problem in Section~\ref{sec:problem}, i.e.,
$\pi(i)=j$ denotes that asteroid $j$ is visited at step $i$ of the
route. Its \emph{ranking}-based counterpart is given by its inverse
$\pi^{-1}$, where $\pi(i)=j \Leftrightarrow \pi^{-1}(j)=i$.  This distinction
is important when calculating distances between permutations, learning a
probability distribution, as done by UMM, or a surrogate model, as done by
CEGO~\citep{IruLop2021gecco}.

The GA used by CEGO to explore the surrogate model uses a budget of $10^4$
evaluations of the surrogate model, population size of 20 individuals, cycle crossover for permutations with crossover rate of 0.5, swap mutation (i.e., exchanging two randomly selected elements) with mutation rate of $1/n$, and tournament selection of size 2 with probability of 0.9. The population used by the GA is different from the solutions evaluated by CEGO because the GA explores the surrogate model, which does not require evaluations of the expensive problem.

We use here the original implementation of
GECO v2.4.2 (\url{https://cran.r-project.org/package=CEGO}). %
UMM is implemented in Python based on the original
code~\citep{IruLop2021gecco}.  We have implemented the proposed ARP benchmark in Python using
the software package \texttt{poliastro} (v0.16)~\citep{poliastro} for
astrophysical calculations.

Each experiment was repeated 30 times with different random seeds.  The
experiments were run in the Picasso supercomputing facility of the University of Málaga with 126 SD530 servers with
Intel Xeon Gold 6230R (26 cores each) at 2.10GHz, 200 GB of RAM and an
InfiniBand HDR100 network.

\subsection{Results of the Black-box Setting}
\label{subsec:black-box}

\begin{table}[!tbp]
  \newcommand{\hdtwo}[1]{\mcol{2}{c}{\bf #1}}
\centering
\caption{Objective value of the best solution found by each algorithm in a single run. Mean (and standard deviation) over 30 independent runs. RS denotes random search.\label{tab:bb_fitness}}
\begin{tabular}{r@{\hskip .75em}*{5}{@{\hskip 1ex plus 1em minus 1em}rr}}
  \toprule
  & \mcol{4}{c}{\bf CEGO}       &\mcol{4}{c}{\bf UMM}\\
  \cmidrule(lr){2-5}\cmidrule(lr){6-9}%
  \bf Instance &\hdtwo{Order} &\hdtwo{Rank} &\hdtwo{Order} &\hdtwo{Rank}         &  \hdtwo{RS} \\
\midrule                                                                                         
  \tt 10\_42 &  379.3& (17.5) &  390.1& (16.4) &  413.6& (15.2) &  388.7& (20.8)  &  419.3& (14.1)\\
  \tt 10\_73 &  346.0& (17.7) &  361.0& (14.3) &  374.5& (16.6) &  361.4& (18.9)  &  375.2& (15.3)\\
  \tt 15\_42 &  575.1& (29.2) &  610.1& (27.4) &  626.9& (26.1) &  594.6& (37.7)  &  636.2& (18.5)\\
  \tt 15\_73 &  582.4& (33.1) &  591.6& (29.1) &  626.2& (24.7) &  595.3& (32.6)  &  636.2& (25.2)\\
  \tt 20\_42 &  806.2& (45.2) &  833.2& (38.7) &  868.9& (26.9) &  841.6& (43.7)  &  877.8& (26.4)\\
  \tt 20\_73 &  839.2& (49.8) &  838.0& (43.0) &  883.0& (28.4) &  876.8& (32.7)  &  897.2& (26.3)\\
  \tt 25\_42 & 1048.8& (51.2) & 1074.5& (56.5) & 1121.2& (28.3) & 1076.9& (48.5)  & 1122.1& (20.8)\\
  \tt 25\_73 & 1069.1& (49.5) & 1096.8& (50.1) & 1131.4& (35.8) & 1115.5& (45.5)  & 1139.7& (25.6)\\
  \tt 30\_42 & 1272.2& (65.8) & 1334.4& (66.6) & 1372.0& (24.3) & 1347.3& (35.6)  & 1357.3& (38.4)\\
  \tt 30\_73 & 1320.9& (63.2) & 1325.1& (66.4) & 1385.9& (27.5) & 1360.6& (51.1)  & 1398.6& (24.9)\\
\bottomrule
\end{tabular}

\end{table}

\begin{figure*}[!tb]
  \centering
  \begin{minipage}{0.67\linewidth}\centering
    \includegraphics[width=\textwidth]{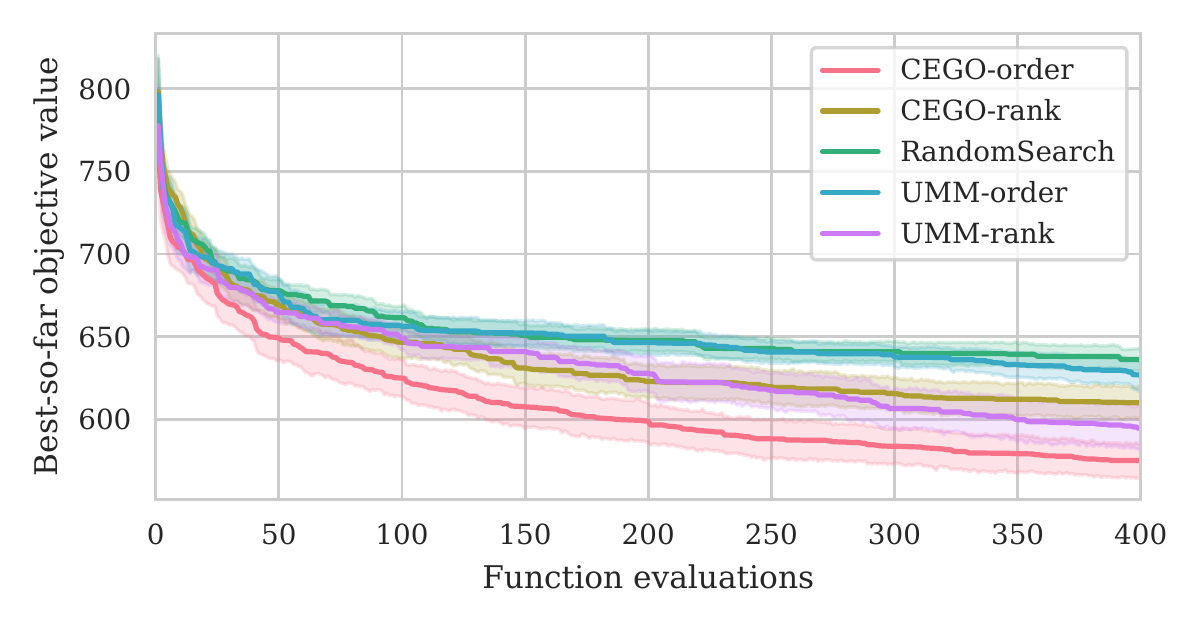}
\end{minipage}
\begin{minipage}{0.67\linewidth}\centering
\includegraphics[width=\textwidth]{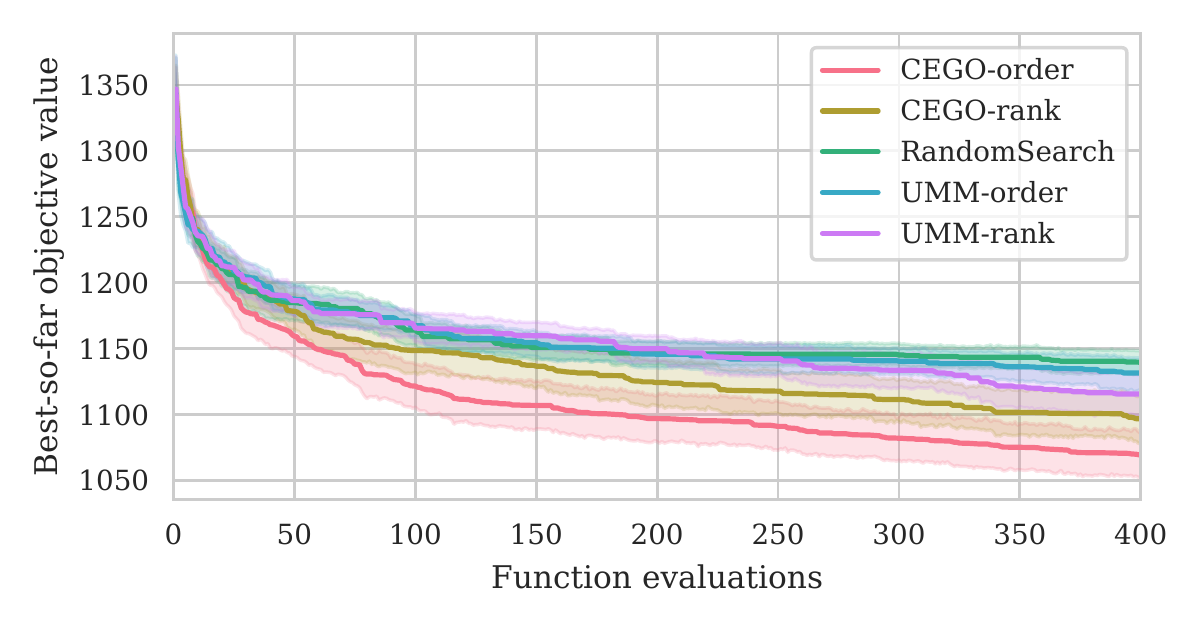}
\end{minipage}
\caption{Convergence of the best-so-far objective value over the number of evaluations. Each line is the mean value over 30 runs and the shaded area is the 95\% confidence interval around the mean. The instances shown are \texttt{15\_42} (top) and \texttt{25\_73} (bottom).}\label{fig:rankVorder}
\end{figure*}

In this section we compare UMM and CEGO (combined with SLSQP) in a black-box
setting to solve the 10 instances of the ARP.  
Table~\ref{tab:bb_fitness} shows the results of CEGO and UMM using the two representations, rank-based or order-based.  As a baseline, we show the results obtained by a simple \emph{random search} that iteratively evaluates 400 random permutations.
We can also see in Figure~\ref{fig:rankVorder} the evolution of the search for two particular instances. Each plot shows as
a solid line the mean, over 30 runs, of the best objective value found up to a
given number of function evaluations. The 95\% confidence interval around the
mean is shown as a shaded region. We applied the non-paired Wilcoxon sumrank test to compare the different algorithm configurations in each instance. The p-values are shown in Table~\ref{tab:pvalues}. We assume that the differences are statistically significant when the p-value is below $0.01$ (confidence level).

Results show that, whenever there are
significant differences (see columns two and three of Table~\ref{tab:pvalues}), rank-based UMM outperforms order-based UMM in terms of the objective function, whereas order-based CEGO outperforms rank-based CEGO. This result matches what we expect given the
assumptions made by each algorithm. That is, UMM internally assumes that the permutation $\pi$ it generates represents a ranking ($\pi(3)=6$ means that object 3 has position~6) and, thus, it makes sense that it works better when the fitness function is aligned with this assumption~\citep{IruLop2021gecco}. CEGO, on the other hand, assumes that the permutation it generates represents an order ($\pi(3)=6$ means that the third object is 6), and that explains the good results of CEGO when the objective function interprets the permutation as an order.

If we compare both algorithms, we observe that CEGO is usually better than UMM using both representations (with only two exceptions for the ranking representation in instances \texttt{10\_42} and \texttt{15\_42}). The best CEGO variant is clearly better than the best UMM variant, although not always with statistical significance (see column four in Table~\ref{tab:pvalues}).

\begin{table}[!bth]
\caption{p-values of a Wilcoxon sumrank test of all the hypotheses checked in during the experimental study for each instance of the problem. Columns two and three correspond to a test comparing the two representations (rank and order) for CEGO and UMM in a black-box setting. The fourth column compares the best black-box UMM (using ranking) with the best black-box CEGO (using order). The final two columns compare UMM-ranking and CEGO-order in a black-box versus informed settings (including the greedy initialization).}
\label{tab:pvalues}
\newcommand{\hls}[1]{\textbf{#1}}
\centering
\begin{tabular}{r@{\hskip 1em}*{2}{@{\hskip 1em plus 1em minus 1em}r}@{\hskip 1em plus 1em minus 1em}c*{2}{@{\hskip 1em plus 1em minus 1em}r}}
  \toprule
             & \mcol{2}{c}{\bf Rank vs. Order} & \mcol{1}{c}{\bf CEGO-order} & \mcol{2}{c}{\bf Greedy vs. BB} \\
  \cmidrule{2-3} \cmidrule{5-6}
\bf Instance & \bf CEGO & \bf UMM & \mcol{1}{c}{\bf vs. UMM-rank} & \bf UMM & \bf CEGO \\ 
  \midrule
  \tt 10\_42 & 6.46e-03 & 2.00e-06 & 4.48e-02 & 2.65e-02 & 1.18e-01 \\ 
  \tt 10\_73 & 2.37e-03 & 1.33e-02 & 7.07e-04 & 8.34e-07 & 1.65e-02 \\ 
  \tt 15\_42 & 3.77e-05 & 6.12e-04 & 5.14e-02 & 2.95e-11 & 1.72e-12 \\ 
  \tt 15\_73 & 3.35e-01 & 2.75e-04 & 1.50e-01 & 6.52e-09 & 3.69e-10 \\ 
  \tt 20\_42 & 1.42e-02 & 8.01e-03 & 2.74e-03 & 2.41e-05 & 3.25e-07 \\ 
  \tt 20\_73 & 7.75e-01 & 4.67e-01 & 2.22e-03 & 2.84e-11 & 2.31e-11 \\ 
  \tt 25\_42 & 3.71e-02 & 7.74e-05 & 3.45e-02 & 3.02e-11 & 4.04e-11 \\ 
  \tt 25\_73 & 4.46e-02 & 1.77e-01 & 5.43e-04 & 2.26e-11 & 2.84e-11 \\ 
  \tt 30\_42 & 5.43e-04 & 4.11e-03 & 1.30e-05 & 2.95e-11 & 3.02e-11 \\ 
  \tt 30\_73 & 6.12e-01 & 3.71e-02 & 1.83e-02 & 2.78e-11 & 3.00e-11 \\ 
   \bottomrule
\end{tabular}
\end{table}

Given these results, we select rank-based UMM and order-based CEGO for the next step in our analysis.

\begin{figure*}[!tb]
  \centering
\begin{minipage}{0.67\linewidth}\centering
\includegraphics[width=\textwidth]{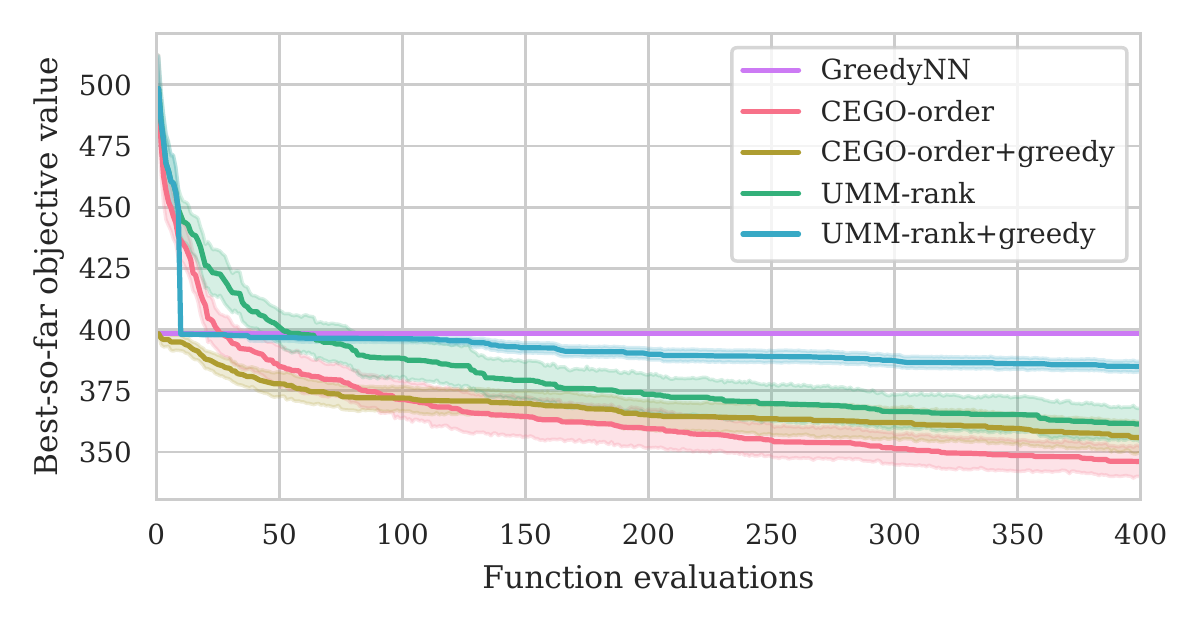}
\end{minipage}
\begin{minipage}{0.67\linewidth}\centering
\includegraphics[width=\textwidth]{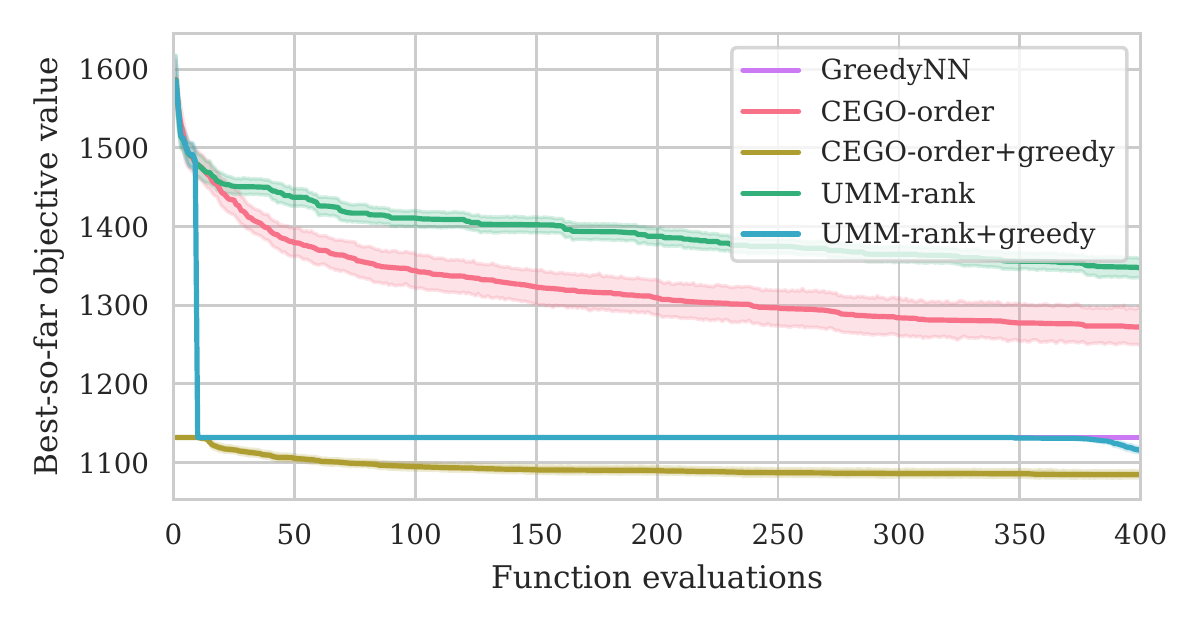}
\end{minipage}
\caption{Convergence of the best-so-far objective value over the number of evaluations. Each line is the mean value over 30 runs and the shaded area is the 95\% confidence interval around the mean. The instances shown are \texttt{10\_73} (top) and \texttt{30\_42} (bottom).}\label{fig:greedy}
\end{figure*}

\subsection{Results of the Informed Setting}
\label{subsec:heuristic}

In this section, we evaluate the effect of providing a good initial solution to
the UMM-rank and CEGO-order variants selected above. This initial solution is
generated using the greedy nearest neighbor heuristic as the initial point of
the max-min-distance sequential design, as explained above, that creates the 10 initial solutions of
both UMM and CEGO.

\begin{table}[!tb]
  \newcommand{\hdtwo}[1]{\mcol{2}{r}{\bf #1}}
\centering
\caption{Objective value of the best solution found by each algorithm in a single run. For CEGO and UMM, the values are the mean (and standard deviation) over 30 independent runs.\label{tab:objvalues}}
\begin{tabular}{r@{\hskip .75em}r*{4}{c@{\hskip 1ex plus 1em minus 1em}rr}}
  \toprule
  &             &\hspace*{.5em}&\mcol{5}{c}{\bf CEGO}       &\hspace*{.5em}&\mcol{5}{c}{\bf UMM}\\
  \cmidrule{4-8}\cmidrule{10-14}
  \bf Inst. &\bf GreedyNN &&\mcol{2}{c}{\bf Order} &\mcol{3}{c}{\bf Ord+Greedy} &&\mcol{2}{c}{\bf Rank} &\mcol{3}{c}{\bf Rank+Greedy} \\
\midrule
  \tt 10\_42 &  391.3 &&  379.3& (17.5) &&      374.9 &(6.9)  &&  388.7 &(20.8) &\hspace*{1em}&     382.4 &(1.1) \\
  \tt 10\_73 &  398.3 &&  346.0& (17.7) &&      355.9 &(18.3) &&  361.4 &(18.9) &&     385.0 &(6.4) \\
  \tt 15\_42 &  508.1 &&  575.1& (29.2) &&      497.2 &(1.1)  &&  594.6 &(37.7) &&     501.7 &(3.8) \\
  \tt 15\_73 &  576.4 &&  582.4& (33.1) &&      525.6 &(4.7)  &&  595.3 &(32.6) &&     545.0 &(9.3) \\
  \tt 20\_42 &  841.7 &&  806.2& (45.2) &&      737.0 &(31.0) &&  841.6 &(43.7) &&     797.6 &(28.3) \\
  \tt 20\_73 &  691.5 &&  839.2& (49.8) &&      661.8 &(7.8)  &&  876.8 &(32.7) &&     684.7 &(6.1) \\
  \tt 25\_42 &  946.3 && 1048.8& (51.2) &&      881.5 &(23.3) && 1076.9 &(48.5) &&     921.7 &(15.2) \\
  \tt 25\_73 &  918.3 && 1069.1& (49.5) &&      873.6 &(15.3) && 1115.5 &(45.5) &&     911.9 &(10.0) \\
  \tt 30\_42 & 1131.7 && 1272.2& (65.8) &&     1084.6 &(13.4) && 1347.3 &(35.6) &&    1115.4 &(11.0) \\
  \tt 30\_73 & 1024.7 && 1320.9& (63.2) &&      967.7 &(17.9) && 1360.6 &(51.1) &&    1017.3 &(6.9) \\
\bottomrule
\end{tabular}
\end{table}

Figure~\ref{fig:greedy} shows the evolution of the algorithms. In
addition to the results provided by CEGO and UMM with and without the initial
greedy solution, we plot the objective function value of the greedy solution as
a horizontal line. 
Table~\ref{tab:objvalues} provides a summary of the results of the two best variants of CEGO and UMM together with the result obtained by the Greedy Nearest Neighbor Heuristic.
A first observation is that UMM only reaches the fitness level of the greedy solution after 10 evaluations. The reason is that we evaluate first the 9 random solutions generated by the max-min-distance
design to show how much the greedy solution improves over them. Moreover, without knowledge of the greedy solution, CEGO or UMM are
only able to match the greedy solution under 400 evaluations on a few small
instances ($n \leq 20$). Unsurprisingly, when starting from the greedy
solution, both algorithms outperform their \emph{uninformed}
counterparts. However, UMM struggles to find any improving solutions, whereas
CEGO is typically able to further improve the greedy solution, although not by
much. In the case of UMM, this behavior is explained by the fact that the
sampling of new permutations is strongly biased by the relative objective
values of the permutations already evaluated. Since the greedy solution is so
much better than any other initial permutation, the algorithm remains stuck at
the initial permutation. In the case of CEGO, starting from a good initial
solution provides a good starting point around which the surrogate model can be
built. However, in several instances, after a quick initial improvement of the
greedy solution, CEGO appears stuck and unable to further find any
improvements. This behavior of CEGO was already reported for other
problems~\citep{IruLop2021gecco} and it is attributed to the inability of the
surrogate model to accurately estimate the underlying objective function beyond
a certain point, thus leading to a blind search.

\begin{table}[!tbp]
\centering
\caption{Mean CPU-time (minutes) of a single run of each algorithm.\label{tab:time}}
\begin{tabular}{r@{\hskip 1em}ccc@{\hskip 1.5ex}ccc@{\hskip 1.5ex}c}
  \toprule
               &             &\hspace*{1em}&\mcol{2}{c}{\bf CEGO}       &\hspace*{1em}&\mcol{2}{c}{\bf UMM}\\\cmidrule{4-5}\cmidrule{7-8}
  \bf Instance &\bf GreedyNN &&\bf Order &\bf Ord+Greedy &&\bf Rank &\bf Rank+Greedy \\
\midrule
   \tt 10\_42 &       0.1 &&      2365.4 &             2147.8 &&      26.0 &             26.2 \\
   \tt 10\_73 &       0.1 &&      2372.9 &             2365.3 &&      24.2 &             25.7 \\
   \tt 15\_42 &       0.1 &&      2445.9 &             2394.1 &&      34.5 &             32.4 \\
   \tt 15\_73 &       0.1 &&      2375.6 &             2404.2 &&      35.5 &             36.2 \\
   \tt 20\_42 &       0.2 &&      2411.4 &             2348.4 &&      47.8 &             48.7 \\
   \tt 20\_73 &       0.2 &&      2223.9 &             2387.4 &&      48.2 &             50.8 \\
   \tt 25\_42 &       0.2 &&      2338.1 &             2452.5 &&      59.4 &             58.9 \\
   \tt 25\_73 &       0.2 &&      2397.3 &             2657.6 &&      47.3 &             60.3 \\
   \tt 30\_42 &       0.2 &&      2375.9 &             2553.7 &&      61.9 &             69.1 \\
   \tt 30\_73 &       0.2 &&      2138.3 &             2611.0 &&      56.0 &             67.4 \\
\bottomrule
\end{tabular}
\end{table}

We also report the mean runtime of CEGO, UMM and the greedy
heuristic in Table~\ref{tab:time}. The runtime of the
greedy algorithm is relatively fast, taking around 6 seconds to generate and
evaluate a single solution. UMM requires between 24 and 69 minutes per run on
average, depending on the instance size, where each run involves 400 solution
evaluations. The runtime per evaluation is larger than in GreedyNN due to the update of and sampling from the probabilistic model.   CEGO also evaluates 400 solutions per run, however, it is more than 100 times slower than UMM, requiring almost 2 CPU-days per run on average for the largest instances tested here.
 Although part of this difference is due to implementation choices, the main difference is the fact that building the Gaussian Process model in CEGO is significantly more expensive than building the probabilistic model in UMM.

\begin{figure}
  \centering
  \begin{minipage}{0.45\linewidth}%
\includegraphics[width=\textwidth]{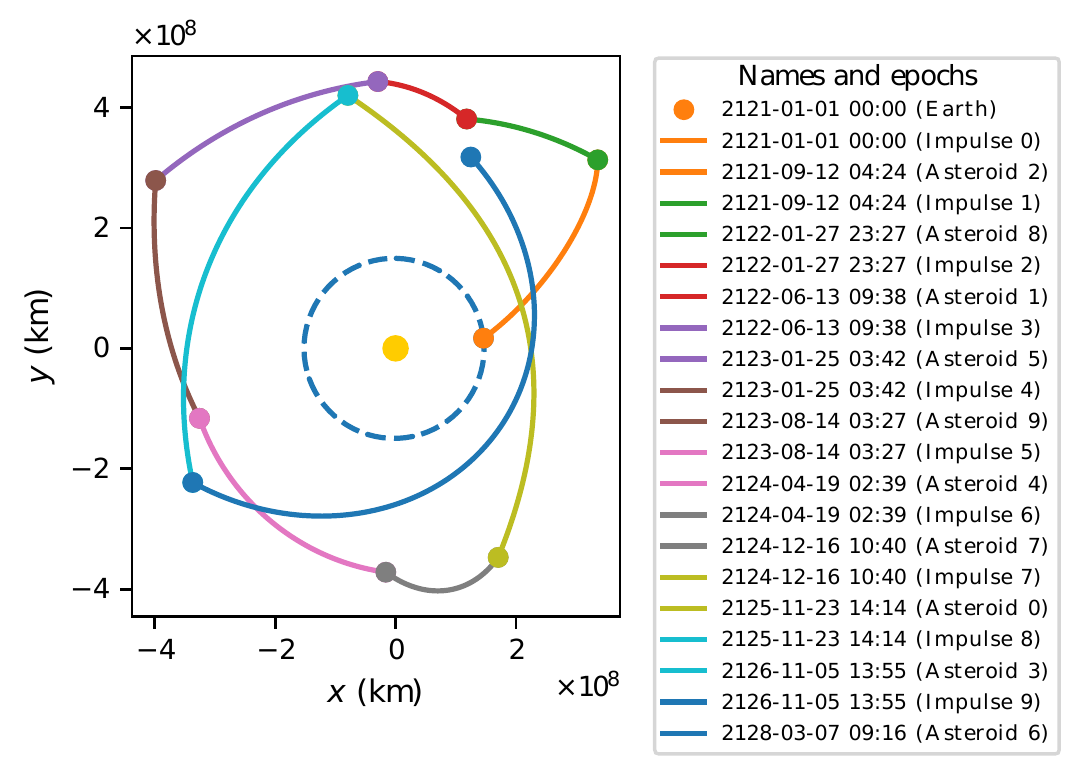}\\[0em]
    {\footnotesize (a)~Greedy Nearest Neighbor.}\hfill
\end{minipage}
  \begin{minipage}{0.45\linewidth}%
    \includegraphics[width=\textwidth]{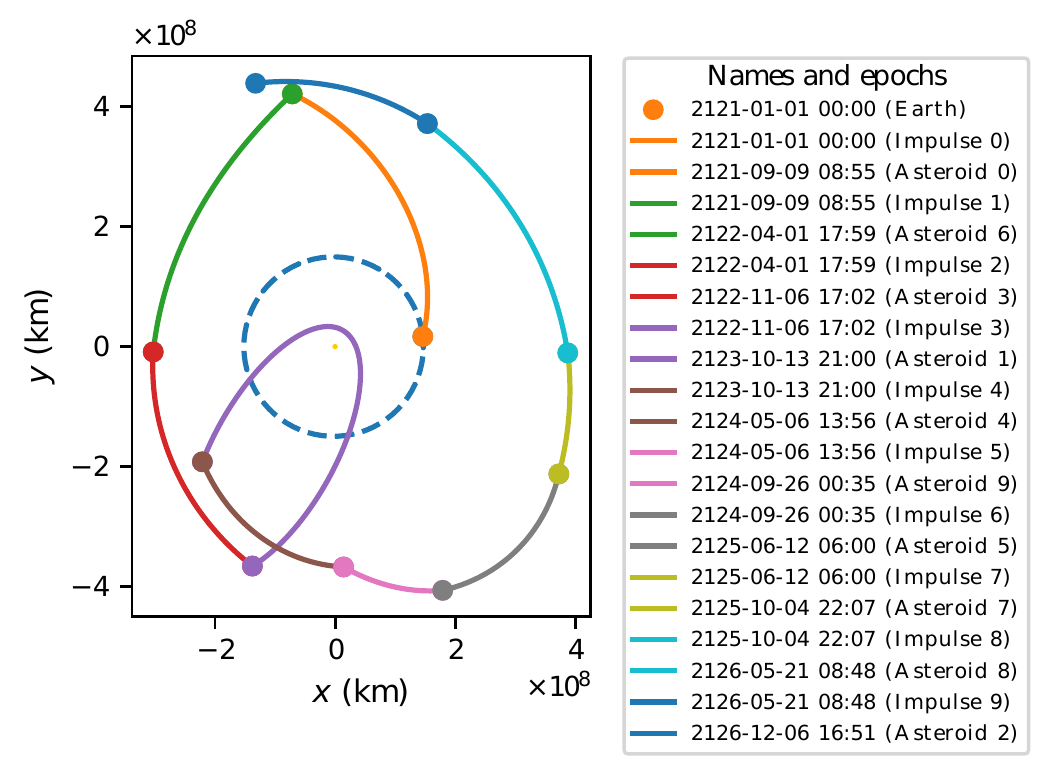}\\[0em]
    {\footnotesize (b) CEGO-order $+$ Greedy.}\hfill
  \end{minipage}
  \begin{minipage}{0.45\linewidth}%
  \includegraphics[width=\textwidth]{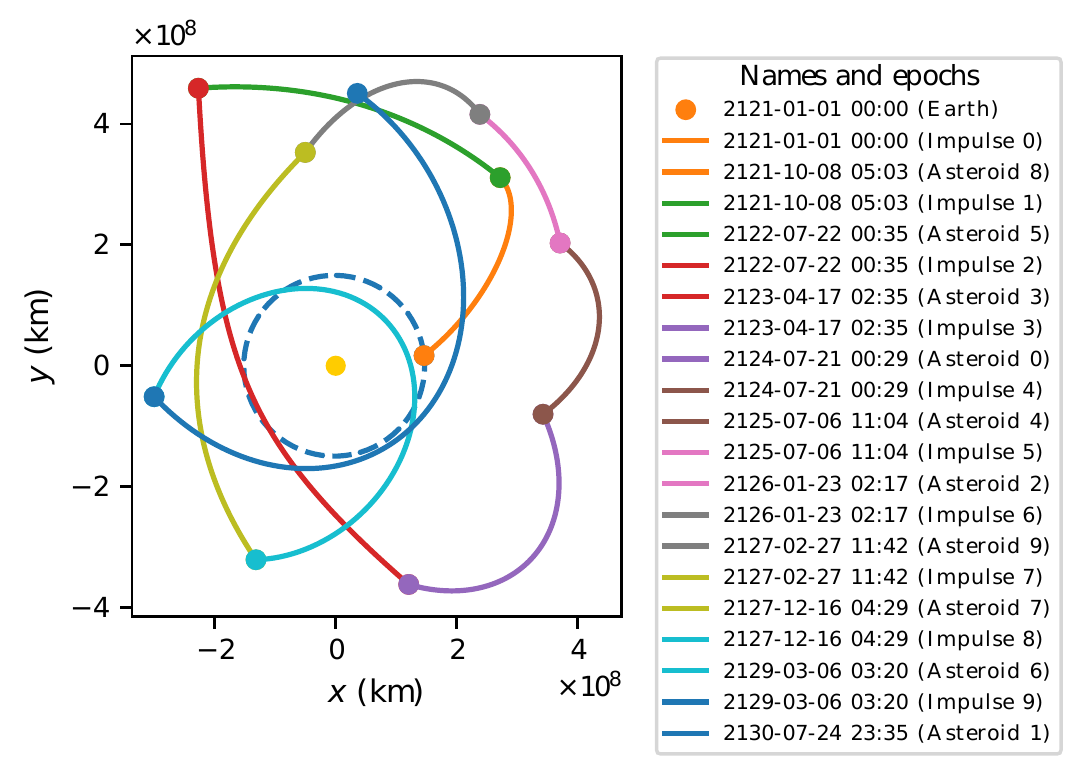}\\[0em]
     {\footnotesize (c) UMM-rank $+$ Greedy.}\hfill
  \end{minipage}
\caption{Example solutions for instance \texttt{10\_73}.}
\label{fig:solutions}
\vspace{-2em}
\end{figure}

Finally, we show in Fig.~\ref{fig:solutions} a visualization of three solutions
found by the greedy heuristic, CEGO-order and UMM-rank, the two last ones
starting from the greedy solution. In each plot, the legend shows the order in
which asteroids are visited (from top to bottom), the epoch at which an impulse
was applied to change orbit, and the arrival epoch at each asteroid. Solid
lines indicate either parking or transfer orbits. When the impulse and arrival
epochs coincide, the parking time at the previous orbit was zero and no parking
orbit is shown. Earth's orbit is shown as a dashed blue line and the Sun's
position as a yellow point. Distances are shown at scale.

\section{Conclusions}\label{sec:conclusion}

We have proposed in this paper a realistic benchmark for evaluating algorithms
designed to tackle expensive black-box optimization problem in permutation
space. 
The results highlight already some limitations of the evaluated
algorithms. First, although CEGO consistently obtains better results than UMM
as long as the right choice of representation is made, this comes at a cost of
100 times more computation time. Second, without any additional information,
very rarely either algorithm is able to match the quality of the solution found
by a greedy approach, suggesting that their results are far from
optimal. Finally, only CEGO makes some use of the information provided by a
good initial solution, whereas UMM is very rarely able to improve over it.

Future work should evaluate other black-box algorithms for permutation problems
in the ARP. Improvements to CEGO and UMM to make better use of good initial
solutions would also be welcome. In fact, it is possible to generate multiple
greedy solutions that may be used as starting points by the black-box
optimizers.  Finally, we plan to change the optimizer of the inner problem from
the current deterministic SLSQP to a stochastic one, which changes the context
of the outer problem from a deterministic permutation problem to a stochastic
one. Under such context, a black-box optimizer for the outer problem must
decide whether it is worth evaluating again already seen permutation in the
hope that a new evaluation may further improve the objective function.

\medskip{}
\begin{smaller}[1]
\noindent\textbf{Reproducibility.} Datasets and source code for reproducing the experiments reported are available from \doi{10.5281/zenodo.5725837}\enspace.\par
\end{smaller}

\medskip{}
\begin{smaller}[2]
\noindent\textbf{Acknowledgements.} This work is partially funded by the Universidad de M\'alaga, Consejer\'ia de Econom\'ia y Conocimiento de la Junta de Andaluc\'ia and FEDER (grant UMA18-FEDERJA-003) and MCIN/AEI/10.13039/501100011033 (grant PID 2020-116727RB-I00). Thanks to the Supercomputing and Bioinnovation Center (SCBI) of the Univerity of Málaga for their provision of computational resources and support. M.\@ L\'opez-Ib\'a\~nez is a ``Beatriz Galindo'' Senior Distinguished Researcher (BEAGAL 18/00053) funded by the Spanish Ministry of Science and Innovation (MICINN).\par
\end{smaller}

\bibliographystyle{splncs04nat}
\providecommand{\MaxMinAntSystem}{{$\cal MAX$--$\cal MIN$} {Ant} {System}}
  \providecommand{\rpackage}[1]{{#1}}
  \providecommand{\softwarepackage}[1]{{#1}}
  \providecommand{\proglang}[1]{{#1}}

\end{document}